\documentclass[10pt]{amsart}
\theoremstyle{plain}
\usepackage{amsmath, amssymb, amscd,graphicx,jrsymb,color,mathrsfs}
\newtheorem{lem}{Lemma}[section]
\newtheorem{prop}[lem]{Proposition}
\newtheorem{thrm}[lem]{Theorem}
\newtheorem{cor}[lem]{Corollary}
\newtheorem{crit}[lem]{Criterion}

\newtheorem*{defn*}{Definition}

\newtheorem{tm}{Theorem}

\theoremstyle{definition}
\newtheorem{defn}[lem]{Definition}

\newcommand{\Z}{{\mathbb{Z}}}
\newcommand{\Q}{{\mathbb{Q}}}

\newcommand{\ts}{{{\thinspace}}}

\newcommand{\hkrn}{HKR_N}
\newcommand{\hkr}{HKR}
\newcommand{\khrn}{\mathscr{P}_N}
\newcommand{\khr}{\mathscr{P}}
\newcommand{\A}{\mathcal{A}}
\newcommand{\thickedge}{\doublepoint}
\newcommand{\wP}{\widetilde{P}}
\newcommand{\wQ}{\mathcal{Q}_N}
\newcommand{\wH}{\widetilde{H}_N}

\newcommand{\wcP}{\widetilde{\mathscr{P}}}
\DeclareMathOperator{\lk}{lk}
\DeclareMathOperator{\Det}{Det}

\title[Khovanov-Rozansky homology of two-bridge knots and links]{Khovanov-Rozansky homology \\ of two-bridge knots and links}
\author{Jacob Rasmussen}
\address{Princeton University Dept. of Mathematics, Princeton, NJ 08544}
\email{jrasmus@math.princeton.edu}
\thanks{The author was supported by an NSF Postdoctoral fellowship.}

\begin{document}

\begin{abstract}We compute the reduced version of
Khovanov and Rozansky's \(sl(N)\) homology for
  two-bridge knots and links. The answer is expressed
  in terms of the HOMFLY polynomial and signature.
\end{abstract}

\maketitle

\section{Introduction}

In \cite{KRI}, Khovanov and Rozansky 
introduced a family of link invariants generalizing the Jones
polynomial homology of \cite{Khovanov}. In this paper, we will mostly
be interested in the reduced version of their homologies, which are
invariants of an oriented
 link \(L\subset S^3\) together with a marked component \(C\)
of \(L\). These invariants take the form of bigraded homology groups
\(\hkrn^{i,j}(L,C)\), where \(N\) is a positive integer. The
information contained in these groups is conveniently represented by
the Poincar{\'e} polynomial
\begin{equation}
\khrn(L,C) =  \sum_{i,j} t^i q^j  \dim \hkrn^{i,j}(L,C).
\end{equation}
Substituting \(t=-1\) gives the graded Euler characteristic, which is
 equal to a classical polynomial invariant of \(L\) --- the \(sl(N)\)
 knot polynomial:
\begin{align}
 \khrn(L,C)|_{t=-1}  = P_L(q^N,q).
\end{align}
Here \(P_L(a,q)\) is the HOMFLY polynomial of \(L\), normalized so
that \(P\) of the unknot is equal to \(1\) and \(P\) satisfies the
skein relation
\begin{equation}
a P(\undercrossing) - a^{-1} P(\overcrossing) = (q-q^{-1})
P(\asmoothing). 
\end{equation}

For \(N=1\), the Khovanov-Rozansky homology is the same for all links
\(L\): \(\khr_1(L,C) \equiv 1\). When \(N=2\), the theory reduces to
the ordinary Khovanov homology, for which extensive computer
calculations have been made by Bar-Natan \cite{DBN}, Bar-Natan and Green
\cite{JavaKh}, and Shumakovitch \cite{KhoHo}.
In contrast,  very few 
computations of the Khovanov-Rozansky homology have been made when
 \(N>2\). Our aim here is to describe the
most elementary type of behavior exhibited by these
theories, and to show that it is satisfied by a  simple class of
links --- the two-bridge links. The result is most easily stated for
knots:
\begin{tm}
\label{Thm:One}
If \(K\) is a two-bridge knot, then for each \(N>4\), 
\(\hkrn(K)\) is
  determined by the HOMFLY polynomial and signature of \(K\). In terms
  of the Poincar{\'e} polynomial, the relation may be expressed as
  follows:
\begin{equation}
\khrn(K) = (-t)^{\sigma(K)/2}P_K(q^Nt^{-1},iqt^{-1/2}).
\end{equation}
\end{tm}

If  \(\khrn(K)\) satisfies the equation above, we say that \(K\)
is \(N\)-thin ({\it cf.} Definition 5.1 of \cite{superpolynomial}.)
Thus the theorem may be summarized by saying that two-bridge knots are
\(N\)-thin for all \(N>4.\) 
The condition that a knot be \(N\)-thin generalizes the definition
of thinness for the ordinary Khovanov homology given in \cite{Khovanov2}. 
Many knots are  known to  be thin in this sense.  In particular, 
 E.S. Lee  proved in \cite{Lee} that any nonsplit alternating
link is thin.  Unfortunately, this fact does not seem to generalize to
the \(sl(N)\) case; in section~\ref{Sec:Examples},
 we give an example of an
alternating knot which is not \(N\)-thin for any \(N > 2\). 
On the other hand, the condition that \(N>4\) is largely 
technical. We expect that two-bridge knots should actually be \(N\)-thin for
all \(N>1\). For \(N=2\), this follows from Lee's theorem, so only the
cases \(N=3,4\) remain unresolved. 

This work in this paper was motivated by the desire to provide some
computational support for the conjectures about the structure of
\(\hkrn\) made in \cite{superpolynomial}. In this regard, we have been
only partly successful, since the knots considered here exhibit
only the simplest possible behavior of the Khovanov-Rozansky
homology. Still, it is worth noting that for the knots
studied in this paper, \(\hkrn\) satisfies both the stabilization
and symmetry properties conjectured in \cite{superpolynomial}. In
addition, we have checked that those knots of \(8\) crossings or fewer
admit plausible candidates for the differentials \(d_1\) and
\(d_{-1}\) described there. 

Although the details of  the proof of Theorem~\ref{Thm:One} are somewhat
messy, the argument itself is quite soft. The main ingredients are the
fact that \(\hkr\) is a link invariant (plus a tiny bit borrowed from
 the actual proof of invariance in \cite{KRI}), the skein exact sequence, and 
the known behavior of \(\hkr\) for the unknot and unlink. 
In particular, 
if the homology groups recently introduced by Khovanov
and Rozansky in \cite{KRII}
could be shown to satisfy a skein exact
sequence for arbitrary diagrams (rather than just braid diagrams), we
expect that an  analogous theorem would hold there  as well. 

The organization of the paper is as follows. The first two sections
contain  background material on singular knots
 and the Khovanov-Rozansky homology.
 In section~\ref{Sec:Thin} we give a more general definition of what
 it means for a knot or link to be \(N\)-thin and describe some
criteria which can be used to prove that a  link is
thin. Finally, in section~\ref{Sec:Examples} we use these criteria
to prove the theorem. We conclude  by
discussing some other knots to which the methods used
in the proof can be applied, and by giving a few calculations of the
unreduced Khovanov-Rozansky homology.  

\noindent{\bf Acknowledgements.} The author would like to thank Nathan
Dunfield, Bojan
Gornik, Sergei Gukov, Mikhail Khovanov, Ciprian Manolescu, Peter
Ozsv{\'a}th and Zolt{\'a}n Szab{\'o} for many helpful discussions
on this subject, and the referee for valuable comments on the
manuscript. 

\section{Invariants of Singular Links}

We begin by fixing our notation for link diagrams and their resolutions.
Suppose \(L\) is an oriented link in \(S^3\) represented by some
planar diagram. Each crossing of the diagram is either {\it
  positive}, like the crossing labeled \(L_+\) in
Figure~\ref{Fig:Diagrams}, or {\it
  negative}, like \(L_{-}\). (Warning: this sign convention is the opposite of
the one used in \cite{KRI}.) We can construct a new oriented link from
\(L\)   either by replacing the crossing with its {\it oriented resolution}
 \(L_0\) or by switching the sign of the crossing.
 If we are willing to forget the orientations, we can also
 form the {\it unoriented resolution} \(L_u\). 

\subsection{Singular links}
\begin{figure}
\includegraphics{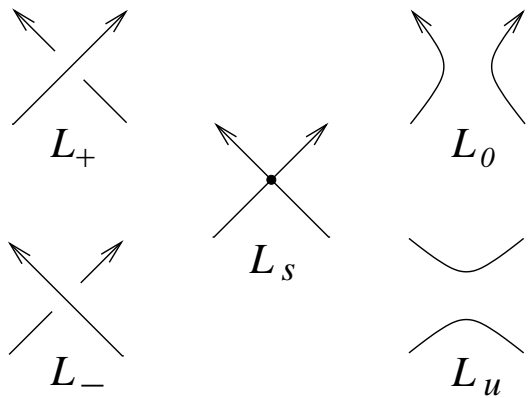}
\caption{\label{Fig:Diagrams} Crossings and resolutions.}
\end{figure}

The original definition of the Khovanov homology involved replacing
each crossing of \(L\) with its oriented and unoriented resolutions. 
To define their more general
 homology theory, Khovanov and Rozansky replaced \(L_u\) with a third
 sort of resolution, in which some 
crossings in the diagram are replaced by the singular diagram
\(L_s\) of Figure~\ref{Fig:Diagrams}.
 In what follows, we will find it convenient to consider the
class of links which contain  one such singular crossing.

\begin{defn}
A {\em singular link} is represented by 
 a planar diagram containing precisely one
 singular crossing.  Two such diagrams
represent the same link if they are  related by a sequence of 
Reidemeister of moves
which take place away from the singular point. 
\end{defn}

In what follows, links represented by diagrams without any singular
crossings will be described as {\it regular}, and the generic term link
will be used to refer to both regular and singular links. If \(L\) is
a link for which we have a specific diagram and crossing in mind,
 we will use \(L_+,L_-,L_s,L_0\) and \(L_u\) to indicate the link in
which the crossing has been replaced by the corresponding crossing or
resolution in Figure~\ref{Fig:Diagrams}. (If we have a
singular link, the crossing to be modified will almost always be the
singular one.)

Many notions from classical knot theory extend naturally to singular
links. We begin with some terminology. 
If \(L\) is a regular link, we let \(c(L)\) denote the number of
components in \(L\). If we modify a crossing of \(L\), we clearly
have \(c(L_-) = c(L_+) = c(L_0) \pm 1\). 
Now suppose \(L\) is singular. The  geometric number of
components of \(L\) is defined to be
\(c(L) = \max\{c(L_-),c(L_0)\}\), and the set of components of
\(L\) is equal to the set of components of the resolution for which the
maximum is attained. 

It is also convenient to consider an invariant \(i(L)\) which reflects
the  algebraic number of
components of \(L\). \(i(L)\) is a \(\Z/2\) valued invariant which
determines things like the parity of the signature and the exponents
of the HOMFLY polynomial. 
If \(L\) is a regular link, it is defined by 
\(i(L) \equiv 1 + c(L) \pmod 2\).
  For singular links, \(i(L)\) is defined to be equal to \(i(L_-)\).

Finally, if \(L\) is a singular link with
components \(C_1\) and \(C_2\), we define the linking number of \(L\) by
\begin{equation}
\lk L = \frac{1}{2}(n_+(C_1,C_2) - n_-(C_1,C_2))
\end{equation}
where \(n_{\pm}(C_1,C_2)\) denotes the number of positive/negative
crossings in \(L\) in which one strand belongs to \(C_1\) and the other to
\(C_2\). Note that \(\lk L \) need not be an integer: if
\(L\) is a singular link, then  \(2  \lk L \equiv i(L) \pmod 2\).

\subsection{The HOMFLY polynomial}

We define the HOMFLY polynomial of a singular link \(L\) by
\begin{align}
\label{Eq:SingHOM}
P(L) & = q P(L_0) - a P(L_-) \\
\label{Eq:SingHOM2}
& = q^{-1}P(L_0) - a^{-1}P(L_+).
\end{align}
(The second equality follows from the skein relation.) Since the link
types of \(L_0,L_-\), and \(L_+\) are invariant under moves that take
place away from the singular point, the HOMFLY polynomial is clearly an
invariant of singular links. 

If \(K\) is a knot,  it is not difficult to see that \(P(K)\)
is a Laurent polynomial in \(a\) and \(q\). In general,
however,   \(P(L)\) is  only a rational function with denominator
\((q-q^{-1})^{c(L)-1}\). 
Most of the links we will consider have two components. 
 In this case the fractional
part is controlled by the following

\begin{lem}
\label{Lem:LinkHom}
 Suppose \(L\) is a two-component link 
 with knots \(K_1\)
  and \(K_2\) as its components. Then
 for any \(n \equiv i(L) \ (\text{\em{mod}} \ 2)\)

\begin{equation}
\wP (L) = P(L) - q^n (-a)^{2 \lk L} P(K_1)P(K_2) \left(
\frac{aq^{-1}-a^{-1}q}{q-q^{-1}}
\right)
\end{equation}
 is a Laurent polynomial in \(a\) and \(q\). 
\end{lem}

\noindent Formulas of this type are well known in the
literature, although we have chosen a somewhat
 nonstandard normalization for the numerator of the
fractional part. 

\begin{proof}
If \(n \equiv m\ (\text{mod} \ 2)\), then \(q^n -q^{m}\) is divisible
by \(q-q^{-1}\). Thus  if the lemma holds for one value of \(n \equiv
i(L) \ (\text{mod} \ 2)\), it
holds for all. Suppose for the moment that \(L\) is a regular
link, so that \(i(L) \equiv 1 \ (\text{mod} \ 2)\).
 If \(L\) is the disjoint union of knots \(K_1\) and \(K_2\), then
\begin{align}
P(L) & = P(K_1)P(K_2)\left( \frac{a-a^{-1}}{q-q^{-1}} \right) \\ & =
P(K_1)P(K_2)\left(a^{-1}q+ q \frac{aq^{-1} - a^{-1}q}{q-q^{-1}}\right)
\end{align}
so the claim holds in this case. Now suppose that \(L_+\) and \(L_-\)
are two diagrams related by a crossing change and that the two strands in
the crossing belong to different components. Then the resolved diagram
 represents a knot \(K\), and \(P(K)\) is a Laurent polynomial. Using
 the skein relation, it is
 easy to see that if one of \(L_-\) or \(L_+\) satisfies the
statement of the lemma, then the other does as well. Since the
components of any link can be unlinked by a sequence of such crossing
changes, the claim holds in the case when \(L\) is a regular link. 

Now suppose that \(L\) is singular, so that one of \(L_0\), \(L_-\) is
a knot and the other is a regular two-component link. Suppose \(L_-\)
is the link. Then \(P(L_0)\) is a Laurent polynomial, so by equation
~\eqref{Eq:SingHOM}, the fractional part of \(P(L)\) is the fractional
part of \(P(L_-)\) multiplied by \(-a\). Clearly \(\lk L  = \lk L_-
+1/2\) and \(i(L) = i(L_-)\), so the lemma holds. On the other hand, 
if \(L_0\) is the link,  the fractional part of \(P(L)\) is the
fractional part of \(P(L_0)\) multiplied by \(q\). Since \(i(L) \equiv
i(L_0) +1 \pmod 2\) and \(\lk L = \lk L_0\), the claim also holds in
this case. 
\end{proof}

\subsection{Determinant and signature}
 The {\it complex
  determinant} of a link \( L\) is defined by
\begin{equation}
\Det L = P_L(-1,i) = V_L(i)
\end{equation} 
where \(V_L(q) = P_L(q^2,q)\) is the Jones polynomial of \(L\). 
If \(\Det L \neq 0 \), we can decompose it into a product of the usual
determinant of \(L\) and a phase \(\phi(L)\):
\begin{align}
\det L & = | \Det L| \\
\phi(L) &= \Det L / |\Det L|.
\end{align}
When it is defined, \(\phi(L)\) amounts to a \(\Z/4\) refinement of
the invariant \(i(L)\), in the sense that \((-1)^{i(L)} =
\phi(L)^2\). 

\begin{lem}
\label{Lem:DetRes}
If \(L\) is a singular link, then \(\det L = \det L_u\).
\end{lem}

\begin{proof}
By equation~\eqref{Eq:SingHOM}, 
\begin{equation}
\Det L = i V_{L_0}(i) + V_{ L_-}(i). 
\end{equation}
On the other hand, Kauffman's unoriented skein relation for the Jones
polynomial tells us that 
\begin{equation}
(-q^{3/2})^{w(L_-)} V_{L_-}(q) = (-q^{3/2})^{w(L_0)} q^{-1/2} V_{L_0}(q) + 
(-q^{3/2})^{w(L_u)} q^{1/2} V_{L_u}(q)
\end{equation}
where \(w(L)\) denotes the writhe of the diagram  \(L\). 
Using the relation \(w(L_0) = w(L_-) +1\), this becomes
\begin{equation}
V_{L_-}(q)+q  V_{L_0}(q) = q^k V_{L_u}(q)
\end{equation}
for some \(k\). 
Substituting \(q=i\) gives \(\Det L  = i^k \Det L_u\), and the claim is
proved. 
\end{proof}

If \(L\) is a regular link, it is well known that the
 phase \(\phi(L)\) can be used to give an inductive characterization
of  the signature 
\(\sigma (L)\) \cite{Conway}.
More precisely, we have

\begin{lem}
\label{Lem:SigPhase}
Suppose \(L\) is a regular link with nonzero determinant. Then 
\(\phi(L) = i^{\sigma(L)}.\)
\end{lem}

\begin{lem}
Let \(L\) be a regular link, and let \(L_0\) be obtained from \(L\) by
resolving a  crossing. If both \(L\) and \(L_0\) have nonzero
determinants, then
 \(\sigma(L_0) = \sigma(L) \pm 1\),
while if one of the two determinants is nonzero and the other is zero,
then 
\(\sigma(L) = \sigma(L_0).\)
\end{lem}
We use this characterization to extend the definition of the signature
to singular links with nonzero determinant. 
\begin{defn}
Let  \(L\) be a singular link with nonzero determinant. Then
 \(\sigma(L)\) is  determined by the requirements that
\(\phi (L)  = i^{\sigma(L)}\)   and 
\(| \sigma(L) - \sigma(L_0)| \leq 1\).
\end{defn}

\begin{cor}
\label{Cor:IndexSig}
If \(\det(L) \neq 0 \), then
 \(i(L) \equiv
\sigma(L) \ (\text{\em{mod}} \ 2)\). 
\end{cor}
\section{Properties of \(\hkr\)}
\label{Sec:KR}
In this section, we briefly review the construction of the
 Khovanov-Rozansky groups
 and describe some of their elementary properties. Our
emphasis is on the formal aspects of the theory --- in particular, we
have suppressed any discussion of matrix factorizations.

\subsection{Foams and functors}
\begin{figure}
\includegraphics{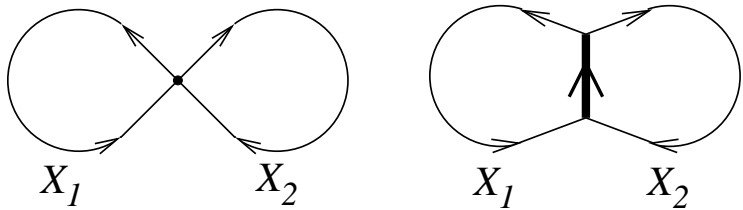}
\caption{\label{Fig:Theta}A simple foam in our notation (left) and
  that of \cite{KRI} (right).}
\end{figure}
The category of \(1\)-manifolds and cobordisms between them plays a
foundational role in the construction of the Khovanov homology. To
define their more general theory,
 Khovanov and Rozansky embedded this category
in a larger one, which we refer to as 
 the category of {\it planar foams}. The objects of this category are
 oriented four-valent planar graphs whose vertices resemble the singular
 crossing of Figure~\ref{Fig:Diagrams}. We allow the possibility that
 some components of the graph are pure circles --- that is, a single
 oriented edge with no vertex. This formulation is slightly different
 from that used \cite{KRI}, where the graphs are trivalent and have
 two kinds of edges. The two are related by the operation of inserting
 a thick edge at each four-valent vertex, as illustrated in
 Figure~\ref{Fig:Theta}. 

\begin{figure}
\includegraphics{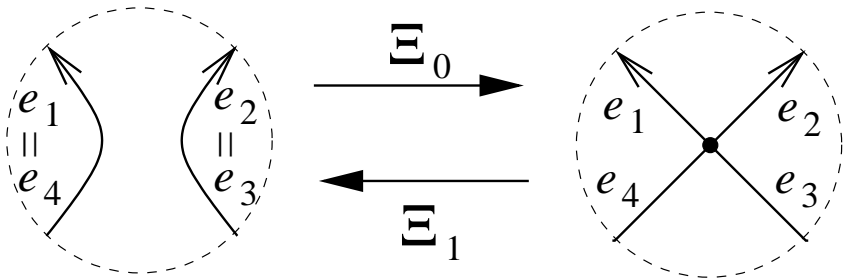}
\caption{\label{Fig:Mor} Elementary morphisms.}
\end{figure}

Morphisms in the category of cobordisms are generated by certain
elementary morphisms (Reidemeister and Morse moves) modulo some
relations (the movie moves of \cite{CarterSaito}.) A similar situation applies
to the category of foams. Rather than describe all the generators and
relations here, we focus on those  which are
relevant for the definition of \(\hkrn\).  For this purpose, 
 it suffices to consider morphisms which are formal compositions
of the elementary morphisms shown in Figure~\ref{Fig:Mor}, 
in which we replace a
region of the graph isomorphic to the region inside one of the dotted
circles with the region inside the other circle. The only relation we
will use is the fact that
morphisms are {\it far-commutative}: if \(\Xi\) and \(\Xi'\) are
morphisms which take place in disjoint circles, then \(\Xi\Xi' =
\Xi'\Xi\). 

 For each integer \(N>0\), Khovanov and Rozansky construct a functor \(\A_N\)
  from  the category of planar foams to
  the category of vector spaces over \(\Q\). These functors play a role
analogous to that of the \(1+1\) dimensional TQFT which appears in the
definition of the usual Khovanov homology. They are defined using the
  theory of matrix factorizations, but we will not  discuss the
details of the construction here. Instead we simply summarize some relevant
facts about \(\A_N\) collected from \cite{KRI}.

 \(\A_N\) is a functor between graded categories. In practical
  terms, this means that \(\mathcal{A}_N(F)\) is a \(\Z\)-graded vector
  space, and  if \(\Xi:F \to
  F'\) is an elementary morphism of the type shown in
  Figure~\ref{Fig:Mor}, then \(\A_N(\Xi) : \A_N(F) \to \A_N(F')\) is a
  graded map of degree one. 
 Another important property of \(\A_N\) is the presence of ``edge
 operators'' associated to the edges of a foam \(F\). More precisely, we have

\begin{prop}
\label{Prop:Foams}
Suppose \(F\) is a foam, and let \(e\) be an
  edge of \(F\). Then there is a linear
operator \(X_e: \A_N(F) \to \A_N(F)\) which satisfies the following
properties:
\begin{enumerate}
\item \(X_e\) is a graded map of degree \(2\). 
\item \(X_e^N=0\). 
\item \(X_eX_{e'} = X_{e'}X_e\) for any two edges  \(e\) and
  \(e'\). 
\item \(\A_N(F)\) is a free module over the ring \(\Q[X_e]/(X^N_e)\). 
\item The action of \(X_e\) commutes with morphisms. In other words,
  if \(\Xi:F \to F'\) is a morphism and \(e\) is a thin edge of \(F\),
  then there is a corresponding edge \(e'\) of \(F'\), and 
\begin{equation*}
\A_N(\Xi) X_e = X_{e'} \A_N(\Xi)
\end{equation*}

\item Let \(e_1, e_2, e_3\), and \(e_4\) be the four  edges adjacent
  to a vertex of \(F\), as illustrated in the right-hand side of
Figure~\ref{Fig:Mor}. Then
\begin{equation*}
X_{e_1}+X_{e_2} = X_{e_3} + X_{e_4} \quad \text{and} \quad X_{e_1}X_{e_2} = X_{e_3}X_{e_4}
\end{equation*}
\item Let \(\A_N(\Xi_0) \) and \(\A_N(\Xi_1) \) be the
  maps associated to the elementary morphisms of
  Figure~\ref{Fig:Mor}. Then
\begin{align*}
\A_N(\Xi_0)  \A_N(\Xi_1) & = X_{e_1} - X_{e_3} \\
\A_N(\Xi_1)  \A_N(\Xi_0)  & = X_{e_1} - X_{e_3}
\end{align*}
where the first operator is viewed as an endomorphism of the foam on
the right-hand side of Figure~\ref{Fig:Mor} 
and the secord is an endomorphism of the foam on the left. 
\end{enumerate}
\end{prop}

We will also need to know the precise structure of \(\A_N(F)\) for some
simple foams.
\begin{lem}
\label{Lem:Theta}
If \(S^1\) is the foam consisting of a single circle with edge
\(e\), then we have \(\A_N(S^1) = \Q[X_e]/(X_e^N)\), and the grading
of \(1\in \A_N(S^1)\) is \(-N+1\).
 If \(\Theta\) is the foam
shown in Figure~\ref{Fig:Theta}, then \(\A_N(\Theta) =
\Q[X_{e_1},X_{e_2}]/(X_{e_1}^N,X_{e_2}^N,S)\), where
\begin{equation}
S = X_1^{N-1}+X_1^{N-2}X_2 + \ldots + X_1X_2^{N-2} + X_2^{N-1}.
\end{equation}
The grading of \(1\in A_N(\Theta)\) is \(-2N+3\). 
\end{lem}

\subsection{The Khovanov chain complex}
\begin{figure}
\includegraphics{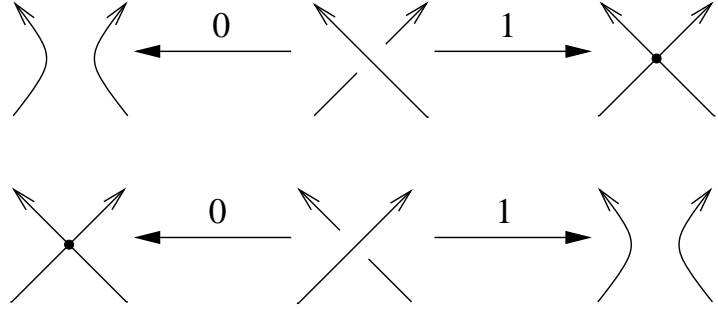}
\caption{\label{Fig:Resolutions} Resolutions of positive and negative
  crossings.}
\end{figure}

Suppose we are given a planar diagram \(D\) representing an
oriented link \(L\), and let \(n\) be the number of nonsingular
 crossings in \(D\).
 The {\it cube of resolutions} of \(D\) is an \(n\)-dimensional cube
 whose vertices are decorated with planar foams and whose edges are
 decorated with morphisms between them. More precisely, each
 nonsingular crossing
 of \(D\) may be resolved in the two ways illustrated in
 Figure~\ref{Fig:Resolutions}. 
If we resolve every crossing of \(D\) in one of these  ways, the
result is clearly a planar foam. Fix once and for all an ordering of
the crossings of \(D\). Then given a vertex \(v\) of the cube
\([0,1]^n\), we associate to it the foam \(F_v\) obtained by resolving
the \(i\)th crossing of \(D\) in accordance with the \(i\)th coordinate
of \(v\). 

Let \(E\) be an edge of the cube with endpoints \(v_0\) and \(v_1\),
where \(v_0\) has one more \(0\) in its coordinates than \(v_1\). We
orient \(E\) so that it points from \(v_0\) to \(v_1\), and write
\(E:v_0 \to v_1\) to indicate this fact. The foams \(F_{v_0}\) and
\(F_{v_1}\) differ only in a neighborhood of a single crossing, where
they resemble the two regions of Figure~\ref{Fig:Mor}. Thus there is
an elementary morphism \(\Xi_E: F_{v_0} \to F_{v_1}\). 

By applying the functor \(\A_N\), we obtain a cube whose vertices are
decorated with graded vector spaces and whose edges are decorated with linear
maps between them. We construct a bigraded chain complex
\(C^{i,j}_N(D)\) from this
cube using the method of \cite{Khovanov}. As a group,  
\(C_N(D) = \oplus_v \A_N(F_v)\). For \(x \in \A_N(F_v)\), the
differential is given by 
\begin{equation}
dx = \sum_{E:v \to v'} (-1)^{s_E} \A_N(\Xi_E)(x)
\end{equation}
where the signs \((-1)^{s_E}\) are chosen so that every two-dimensional
face of the cube has an odd number of minus signs. (This ensures that
\(d^2=0\).) There are many ways to do this, but they all result in
isomorphic chain complexes. 

The bigrading on \(C^{i,j}_N(D)\) is defined as follows. Let \(x\) be
a homogenous element of 
\(\A_N(F_v)\) with grading \(q(x)\). Then 
\begin{align}
\label{Eq:iGr}
i(x) & = s(v)-n_+(D) \\
\label{Eq:jGr}
j(x) &= q(x) -i(x) +(N-1)(n_+(D) - n_-(D))
\end{align}
 where \(s(v)\) denotes the sum of the coordinates of \(v\) and
 \(n_+\) and \(n_-(D)\) are the
number of positive and negative crossings in \(D\). The gradings are
chosen so that \(i(dx) = i(x) + 1\) and \(j(dx) = j(x)\). 

Khovanov and Rozansky showed that the graded Euler characteristic of 
\(C_N^{i,j}(D)\) is the unnormalized \(sl(N)\) polynomial of the link
\(L\) it represents. In
other words, 
\begin{equation}
\sum_{i,j} (-1)^i q^j \dim C_N^{i,j}(D) = \frac{q^N-q^{-N}}{q-q^{-1}}
P_L(q^N,q) 
\end{equation}
(To be strict, the result of \cite{KRI} is fomulated for ordinary
links. Our definition of the HOMFLY polynomial of a singular link was
chosen so that the formula holds in this case as well.) 

In contrast to the case of the ordinary Khovanov homology, where the
fact that the graded Euler characteristic gives the Jones polynomial
is almost immediate from  the construction, the proof of this theorem
is far from trivial.  The argument uses a familiar state model for the
\(sl(N)\) polynomial \cite{MuOhYa} (see also \cite{JonesHOM}
\cite{TuraevHOM}) but also requires a careful analysis of the
properties of the functor \(\mathcal{A}_N\). 

The second major result of \cite{KRI} is that the homology of the
chain complex \(C_N^{i,j}(D)\)
depends only on \(L\), and not on the particular
planar diagram \(D\) that we used to represent it. This group is the
unreduced Khovanov-Rozansky homology of \(L\), and is denoted by
\(H_N^{i,j}(L)\). The theorem is proved by checking that the homology
is invariant under the Reidemeister moves, so it applies equally well
to both regular and singular links. 

\subsection{Skein Exact Sequences}

Let \(D\) be a planar diagram of a regular 
 link \(L\), and let \(c\) be a crossing  of \(D\). 
 Then it is easy to see that there is a short exact sequence
\begin{equation*}
\begin{CD}
0 @>>> C_N(D_1) @>>> C_N(D) @>>> C_N(D_0) @>>> 0 
\end{CD}
\end{equation*}
where \(D_0\) is the planar diagram in which \(c\) has been
given the \(0\) resolution, and \(D_1\) is the diagram in which \(c\)
has been given the \(1\)-resolution. 
We call the resulting long exact sequence on homology a {\it skein
  exact sequence.}  These sequences are the main computational tool
used in the proof of Theorem~\ref{Thm:One}.
For future reference, we  record 
the degrees of the maps involved in them. 

\begin{lem}
\label{Lem:Skein}
There are long exact sequences
\begin{equation*}
\begin{CD}
\cdots @>>>H_N(L_s)@>{(1,-N)}>>H_N(L_-) @>{(0,N-1)}>> 
H_N(L_0) @>{(0,1)}>>
H_N(L_s) @>>> \cdots
\end{CD}
\end{equation*}
and
\begin{equation*}
\begin{CD}
\cdots @>>> H_N(L_0) @>{(0,N-1)}>> H_N(L_+) @>{(1,-N)}>>
H_N(L_s) @>{(0,1)}>> H_N(L_0) @>>> \cdots
\end{CD}
\end{equation*}
\end{lem}
The numbers in parentheses over each map indicate its effect on the
bigrading. For example, if \(x \in H_N(L_s)\), then the first
map in the first sequence takes \(x\) to an element \(y \in
H_N(L_-)\), and 
\begin{equation}
i(y)  = i(x) + 1 \quad \quad
j(y)  = j(x) - N. 
\end{equation}
The proof of the lemma is a straightforward computation using
equations~\eqref{Eq:iGr} and \eqref{Eq:jGr}. 

\subsection{Reduced homology}
\label{SubSec:Reduced}

Let \(e\) be an edge of the planar diagram \(D\), and 
let \(F_v\) be a complete resolution of
\(D\). Then \(e\) determines an  edge \(e_v\) of \(F_v\), and 
we give the group \(C_N(D)\) the structure of a module over \(\A[X_e]\) by
defining \(X_e(x) = X_{e_v}(x)\) for \(x \in \A_N(F_v)\). 
If \(E:v \to v'\) is an edge of the cube of resolutions, part (5) of
Proposition~\ref{Prop:Foams} implies  that \(\mathcal{A}_N 
(\Xi_E) X_{e_v} = X_{e_{v'}}
\mathcal{A}_N (\Xi_E)\). Thus \(d\) commutes with the action of \(X_e\), and the chain complex
\(C_N(D)\) is a module over \(\Q[X_e]\). 

The reduced homology \(\hkrn(L,C)\) discussed in the introduction is
defined to be the homology of the chain complex
\begin{equation}
C_N(D)/(X_eC_N(D)) \cong C_N(D) \otimes_{\Q[X_e]} \Q
\end{equation}
where \(e\) is any edge of \(D\) belonging to \(C\). The homology is
independent of the choice of \(e\), since any diagram of \(L\) with a
marked point on \(C\) can be transformed into any other 
by a sequence of Reidemeister moves and isotopies in
\(S^2\) which take place away from the marked point. 

Applying the decomposition of the previous section to the reduced chain
complex, we see that the reduced homology satisfies skein exact
sequences analogous to those of Lemma~\ref{Lem:Skein}. Note that
these sequences depend not only on the choice of diagram \(D\), but
also on the choice of the edge \(e\) in \(D\) with respect to which we
reduce. Varying \(e\) will result in exact sequences involving the
reduced homology with respect to different components. 

 The module structure of \(C_N(D)\) makes \(H_N(L)\) into a module
 over \(\Q[X_e]\). 

\begin{lem}
Suppose \(e\) and \(e'\) belong to the same component of \(L\). Then
their action on \(H_N(L)\) is the same. 
\end{lem}

\begin{proof}
It suffices to prove the result when \(e\) and \(e'\) are two edges
separated by a single crossing. We assume that the edges near the
crossing are labeled as in the right-hand side of 
Figure~\ref{Fig:Mor}. 

Let \(a\)
be a closed element of \(C_N(D)\), and write \(a = b+c\) for \(b
\in C_N(D_0) \) and \(c \in C_N(D_1) \), where \(D_0\) and \(D_1\) are
the diagrams obtained by giving the \(0\) and \(1\) resolution to the
crossing in question. Then we can write
\begin{equation}
d a  = d_0 b  + d_{01}b + d_1c
\end{equation}
where \(d_i\) is the differential in \(C_N(D_i)\) and \(d_{01}:
C_N(D_0) \to C_N(D_1)\). Since \(a\) is closed, it follows that \(d_0b
= 0 \) and  \(d_{01}b = -d_1c\). 

The action of the map \(d_{01}\) may be described as follows.
 Suppose \(v_0\) is a vertex of the cube of resolutions of
 \(D_0\). Then there is a corresponding vertex \(v_1 \) in the cube of
 resolutions of \(D_1\), and a unique edge \(E:v_0 \to v_1\). If \(x
 \in \A_N(F_{v_{0}})\), then 
\begin{equation}
d_{01}x = (-1)^{s_E} \A_N(\Xi_E)(x).
\end{equation}
Let \(\overline{E}:v_1\to v_0\) be the same edge with the opposite
orientation, and consider the map \(d_{10}:C_N(D_1) \to C_N(D_0) \)
given by 
\begin{equation}
d_{10}(y) = (-1)^{s_E} \A_N(\Xi_{\overline{E}})(y).
\end{equation}
Then 
\begin{equation}
d_{10}d_{01}x = \A_N(\Xi_{\overline{E}})\A_N(\Xi_E)(x) =
(X_{e_1}-X_{e_3})x
\end{equation}
by property (7) of Proposition~\ref{Prop:Foams}. Similarly, we find
that \(d_{01}d_{10}y = (X_{e_1}-X_{e_3})y\) for \(y\in
C_N(D_1)\).

Consider  \(d_{10}c\) as an element of \(C_N(D)\). We compute
\begin{align}
dd_{10}c &= d_0d_{10}c + d_{01}d_{10}c \\
& = -d_{10}d_1c + d_{01}d_{10}c \\
& = d_{10}d_{01}b + d_{01}d_{10}c \\
& = (X_{e_1}-X_{e_3})b + (X_{e_1}-X_{e_3})c \\
& = (X_{e_1}-X_{e_3})a. 
\end{align}
Thus \(X_{e_1}a \) is homologous to \(X_{e_3}a\). If the pair of
edges in question was \(e_1\) and
\(e_3\), the claim is proved. On the other hand, if the pair of edges
was \(e_2\) and \(e_4\), we apply property (6) of
Proposition~\ref{Prop:Foams} to get
\(
(X_{e_1}-X_{e_3})a = (X_{e_4}-X_{e_2})a.
\)
We then argue as before. 
\end{proof}

It follows that \(H_N(L)\) can be naturally viewed as
 a module over \(\Q[X_i]\), where
\(i\) runs over the set of components of \(L\). The proof of 
 the lemma carries over verbatim to the reduced complex
 \(C_N(D)/(X_eC_N(D))\), so a similar result holds for the reduced
 homology. In this case, however, the action of the edge \(X_e\)
 is tautologically \(0\), so we have
\begin{cor}
\(X_e\) acts by \(0\) on \(\hkrn(L,C)\) whenever \(e\) belongs to \(C\).
\end{cor}

\section{Thin Knots and Links}
\label{Sec:Thin}

In this section, we define what it means for a knot or a two-component
link (regular or singular) to be thin, and describe the basic
properties of such links. We then give some criteria which can be used
to recognize thin links.  

\subsection{Thin knots}

We have already given one formulation of what it means for a knot to be
\(N\)-thin in the introduction. The definition we give here is less compact 
but perhaps more illuminating. 

\begin{defn}
\label{Def:Alt}
Let \(P(a,q) = \sum c_{mn} a^mq^n\) be a Laurent polynomial in
\(a\) and \(q\). We
say that \(P\) is {\it alternating} if each nonzero term in the sum
 \(P(-1,i) = \sum c_{mn} (-1)^m i^n \)
has the same phase. 
\end{defn}

Suppose that \(K\) is a knot and that \(P(K) = \sum c_{mn}
a^mq^n\) is its HOMFLY polynomial. Then \(c_{mn} = 0\) unless \(m\)
and \(n\) are both even, so the condition that \(P(K)\) be alternating
amounts to saying that the sign of \(c_{mn}\) is determined by the
parity of \(n/2\). Starting from \(P(K)\), we
 form the three-variable polynomial
\begin{equation}
\label{Eq:Thin}
\mathscr{P}(K) = \sum |c_{mn}| \ts a^mq^n t^{(\sigma(K) - 2m - n)/2}.
\end{equation}
The coefficients of \(\mathscr{P}(K)\) are all positive, so it is
potentially a Poincar{\'e} polynomial. It is not difficult to see that
if we substitute \(t=-1\) in \(\mathscr{P}(K)\), we recover \(P(K)\)
if and only if \(P(K)\) is alternating.

\begin{defn} 
\label{Def:KnotThin}
For \(N>2\), a knot \(K\) is \(N\)-thin if \(P(K)\) is alternating and
\begin{equation}
\khrn(K) = \mathscr{P}(K)|_{a=q^N}. 
\end{equation}
\end{defn}

\begin{cor}
If \(K\) is \(N\)-thin, then \(\dim \hkrn(K) = \det K\). 
\end{cor}

\begin{proof}
\(P(K)\) is alternating, so
\begin{equation}
\det K  = |P_K(-1,i)| 
 = \biggl\vert \sum c_{mn}(-1)^mi^n \biggr \vert 
 = \sum |c_{mn}|.
\end{equation}
since all terms in the first sum have the same phase. 
\end{proof}

{\bf \noindent Remark:} If \(K\) is thin, then
up to a change of variables (due to differing choices of
normalization for \(\hkrn\)) \(\mathscr{P}(K)\) is
the superpolynomial described in  \cite{superpolynomial}.

Of course, we should check that the definition given above
agrees with the one used in the introduction. 

\begin{lem}
 For \(N>2\), \(K\) is \(N\)-thin if and only if 
\begin{equation}
\label{Eq:ThinDef}
\khrn(K) = (-t)^{\sigma(K)/2} P_K(q^Nt^{-1},iqt^{-1/2}).
\end{equation}
\end{lem}

\begin{proof}
Let
 \(P(K) = \sum c_{mn}a^mq^n\). If \(P(K)\) is alternating, then
\begin{equation*}
(-t)^{\sigma(K)/2} P_K(at^{-1},iqt^{-1/2})   = (-1)^{\sigma(K)/2}\sum
  (-1)^{n/2}c_{mn} a^mq^nt^{(\sigma(K)-2m-n)/2} 
 = \mathscr{P}(K)
\end{equation*}
since the sign of \((-1)^{n/2}c_{mn}\) is given by
\(\phi(K)\), and  by Lemma~\ref{Lem:SigPhase}  this is also equal to
\((-1)^{\sigma(K)/2}\).  
If \(K\) is \(N\)-thin, then \(P(K)\) is alternating by definition, and
substituting \(a\)=\(q^N\) gives the desired result. 
Conversely, suppose  equation~\eqref{Eq:ThinDef} holds.
 Then a term \(c_{mn}a^mq^n\)
in \(P(K)\) gives rise to a term 
\begin{equation}
C_{mn}=(-1)^{(\sigma(K)+n)/2} c_{mn} q^{Nm+n}t^{(\sigma(K) - 2m -n)/2}
\end{equation}
 in \(\khrn(K)\). For \(N\neq 2\), different values of \(m\)
 and \(n\) always give rise to different exponents of \(q\) and \(t\)
 in \(C_{mn}\),
 so the sign of each individual term \(C_{mn}\) must be positive. 
 Thus \((-1)^{\sigma(K)+n/2}c_{mn} \geq 0\), which implies that 
 \(P(K)\) is alternating and \(\khrn(K) = \mathscr{P}(K)|_{a=q^N}. \)
\end{proof}

When \(N=2\), the two formulations diverge. We leave it to the
reader to check that in this case, equation~\eqref{Eq:ThinDef} is
satisfied if and only if the usual Khovanov homology is thin in the
sense of \cite{Khovanov2} and the invariant \(s(K)\) described  in 
\cite{khg} is
equal to \(\sigma(K)\). In contrast, Definition~\ref{Def:KnotThin}
imposes the additional constraint that the HOMFLY polynomial of \(K\)
be alternating. An example of a knot satisfying the first two conditions
but not the third is given in section~\ref{Sec:Other}.

Another useful characterization of thinness is in terms of the
\(\delta\)-grading defined in
 \cite{superpolynomial}. If \(c_{lmn}t^la^mq^n\)
is a monomial in \(\mathscr{P}(K)\), we assign to it the
\(\delta\)-grading \(\delta = 2l+2m+n\). Then equation~\eqref{Eq:Thin}
may be described by saying that if \(K\) is thin, all terms
in \(\mathscr{P}(K)\) have \(\delta = \sigma(K)\). 
The corresponding grading on \(\hkrn\) is 
defined by \(\delta(x) = 2i+j\) for \(x \in \hkrn^{i,j}(K)\). Although
\(\delta(x)\) is well-defined as an integer,
 it is best viewed as an element of \(\Z/(N-2)\Z\). Indeed,
 substituting \(a=q^N\) turns the monomial \(t^l a^m q^n\), which has 
 \(\delta\)-grading \(2l+2m+n\)) into \(t^l q^{n+mN}\), which has 
\(\delta\)-grading \(2l+n+mN\), and the two quantities agree modulo
\(N-2\).  Thus if \(K\) is \(N\)-thin,   \(\delta(x) \equiv \sigma(K)
\pmod{N-2}\) for all \(x \in \hkrn(K)\).

\subsection{Thin links}
For a number of reasons, the definition of thinness for links is more
complicated than for knots. 
First, there is the question of the module structure.
 Recall from section~\ref{SubSec:Reduced} that \(\hkrn(L,C)\) is a
 module over \(\Q[X_i]\), where \(i\) runs over the set of components
 of \(L\). The variable corresponding to \(C\)  acts by
\(0\), so in the case of a knot this issue does not arise. 
To describe the homology associated to a thin
 two-component  link, however, we must specify not only its Poincar{\'e}
 polynomial, but also its structure as a \(\Q[X]\) module, where
 \(X\) is the variable corresponding to the unmarked component of
 \(L\).

Second, if \(L\) is a link with more than one component, \(P(L)\) is not a
Laurent polynomial, and Definition~\ref{Def:KnotThin} cannot be
applied. To simplify matters, we assume that 
 \(L\) is a two-component link (so its HOMFLY polynomial is controlled
 by  Lemma~\ref{Lem:LinkHom}) and that both components of \(L\) are
 unknots. If we further suppose that \(\det L \neq 0 \), we can apply
 Lemma~\ref{Lem:LinkHom} and Corollary~\ref{Cor:IndexSig} to write
\(P(L) = \wP(L) + Q(L)\), where 
\begin{equation}
\label{Eq:ThinLink}
Q(L) =  q^{\sigma(L)}(-a)^{2 \lk L} \left( \frac{aq^{-1}
  -a^{-1}q}{q-q^{-1}} \right)
\end{equation}
and \(\wP(L)\) is a Laurent polynomial in \(a\) and \(q\). 
Write \(\wP(L) = \sum c_{mn}a^mq^n\). In analogy with
equation~\eqref{Eq:Thin}, we set
\begin{equation}
\widetilde{\mathscr{P}}(K) = \sum |c_{mn}| \ts a^mq^n t^{(\sigma(L) - 2m - n)/2}.
\end{equation}

\begin{defn} Let \(L\) be a two-component link  both of whose
  components are unknots and with \(\det L \neq 0 \). 
 If \(C\) is a component of \(L\),
 we say that the pair \((L,C)\) is \(N\)-thin if
\begin{enumerate}
\item \(\hkrn(L,C) \cong \wH(L,C) \oplus \Q[X]/(X^{N-1})\), where the action
  of \(X\) on \(\wH\) is trivial.
\item \(\wP(L)\) is alternating and the Poincar{\'e} polynomial of
  \(\wH(L,C)\) satisfies
\begin{equation*}
\mathscr{P}(\wH(L,C)) = \widetilde{\mathscr{P}}(L)|_{a=q^N}.
\end{equation*}
\item The Poincar{\'e} polynomial of the second summand is given by
\begin{equation*}
\wQ(L) = 
q^{\sigma(L)} (q^Nt^{-1})^{2 \lk L}(q^{-N+2}+q^{-N+4} + \ldots
q^{N-2}).
\end{equation*}
\end{enumerate}
\end{defn}

\begin{cor}
If \((L,C)\) is \(N\)-thin, then \(\dim \hkrn(L,C) = \det L +N -2\). 
\end{cor}

\begin{proof}
If we substitute \(a=-1\) and \(q=i\) in \eqref{Eq:ThinLink}, 
the second term reduces to \(i^{\sigma(L)} \), which has the same
phase as \(\Det L\). It follows that 
\begin{equation}
\det L  = 1 + |\wP_L(1,i)| 
       = 1 + \sum |c_{mn}| 
       = 1 + \dim \wH(L,C).
\end{equation}
The second summand clearly has dimension 
\(N-1\), so the claim follows. 
\end{proof}

If \((L,C)\) is \(N\)-thin, then every  element in the first summand
has \(\delta\)-grading \(\sigma(L)\), just as it is for knots.
 The \(\delta\)-grading
of the generators of
 the second summand varies, but it is easy to see that the
 terms with the highest and lowest
\(q\)-gradings also have \(\delta\)-grading congruent to \(\sigma(L)
\pmod{N-2}\).  

\subsection{Exact sequences}
We are now in a position to state our  main technical result. 

\begin{thrm}
\label{Thm:ThinSeq}
Suppose \(L_1\), \(L_2\), and \(L_3\) are a knot, a regular
two-component link, and a singular two-component link (not necessarily
in that order) related by a skein exact sequence, and that 
\begin{equation}
\det L_2 = \det L_1 + \det L_3.
\end{equation}
If \(L_1\) and \(L_3\) are \(N\)-thin for some \(N>4\), then
\(L_2\) is \(N\)-thin as well. 
\end{thrm}

In the interest of maintaining a uniform notation for knots and links,
 we have omitted mention of the marked
 components. For those components which are links, the statement
 should be taken to 
 refer to the reduced homology with respect to the component used to
 define the skein exact sequence. 

\begin{proof}
Without loss of generality, we assume that the sequence is
arranged as follows:
\begin{equation*}
\begin{CD} @>>> \hkrn(L_1) @>{f_1}>> \hkrn(L_2) @>{f_2}>> \hkrn(L_3) @>{f_3}>>
  \hkrn(L_1) @>>>
\end{CD}
\end{equation*}
For \(x \in \hkrn^{i,j}(L_n)\), we define \( \Delta(x) \in
\Z/(N-2)\) by
\begin{align}
\Delta (x) & = \delta(x) -
\sigma(L_n) \\ & = 2i+j- \sigma(L_n).
\end{align}

\begin{lem}
\label{Lem:Delta}
Under the hypotheses of the theorem, 
\(\Delta\) is preserved by \(f_1\) and \(f_2\), while \(f_3\) raises
\(\Delta\) by \(2\). 
\end{lem}

\begin{proof}
The effect of the maps
 \(f_i\) on the \(\delta\)-grading is easily determined from
 Lemma~\ref{Lem:Skein}. It is given by 
\begin{equation*}
\begin{CD}
\ @>>>\hkrn(\thickedge)@>{0}>>\hkrn(\undercrossing) @>{1}>> 
\hkrn(\asmoothing) @>{1}>>
\hkrn(\thickedge) @>>> 
\end{CD}
\end{equation*}
\begin{equation*}
\begin{CD}
 @>>> \hkrn(\asmoothing) @>{1}>> \hkrn(\overcrossing) @>{0}>>
\hkrn(\thickedge) @>{1}>> \hkrn(\asmoothing) @>>> 
\end{CD}
\end{equation*}
where the number over each arrow indicates the degree by which it
raises \(\delta\). To determine the relation between the signatures,
we use the relations
\begin{align}
\Det(\doublepoint) & = i \Det (\asmoothing) + \Det(\undercrossing) \\
\Det(\doublepoint) & = -i \Det (\asmoothing) + \Det(\overcrossing) 
\end{align}
obtained by substituting \(a=-1\), \(q=i\) in
equations~\eqref{Eq:SingHOM} and \eqref{Eq:SingHOM2}.
 For example, if \(L_1 = \doublepoint\),
\(L_2 = \undercrossing\), and \(L_3=\asmoothing\), then in order to have
 \(\det \undercrossing = \det \doublepoint + \det \asmoothing\), we
 must have \(\phi(\undercrossing) = - i \phi(\asmoothing) =
 \phi(\doublepoint)\). Thus \(\sigma(\undercrossing) =
 \sigma(\asmoothing)-1 = \sigma(\doublepoint)\), and a quick
 comparison with the first exact sequence above verifies the claim of
 the lemma. We leave it to the reader to check the remaining five
 cases, which are all similar. 
\end{proof}

By hypothesis, \(L_1\) and \(L_3\) are \(N\)-thin, so we can write
\(\hkrn(L_n) \cong A_n \oplus B_n\) (\(n=1,3\)), 
where every element of \(A_n\) has \(\Delta\)-grading
\(0\) and \(B_n \) is trivial if \(L_n\) is a knot and isomorphic to 
\(\Q[X]/(X^{N-1})\) if \(L_n\) is a link.

\begin{lem}
 \(f_3=0\) unless \(L_2\) is a knot, in which case \(f_3\) acts
 trivially on \(A_3\) and sends 
 \(B_3 \)  to 
 \(B_1 \) by multiplication by \(cX\) for some
 \(c\neq 0 \). 
\end{lem}

\begin{proof}
We consider the various
 components of \(f_3\) with respect to the direct sum decompositions
 of \(\hkrn(L_1)\) and \(\hkrn(L_3)\). 
We start with the component which maps \(A_3\) to \(A_1\). By
the provious lemma, we know that \(\Delta(f_3(A_3)) \equiv 2 \pmod
{N-2}\), while \(\Delta(A_1) \equiv 0\). Thus for \(N>4\), this
component must be trivial. 

Next, consider the component mapping
\(A_3\) to \(B_1\) (if it exists). Since \(f_3\) is a map of \(\Q[X]\)
modules and \(X\) acts trivially on \(A_3\), the image of this map
must be spanned by \(X^{N-2} \in B_1 \cong  \Q[X]/(X^{N-1})\). But
\(X^{N-2}\) also has \(\Delta\)-grading congruent to \(0\), so this
component must be trivial as well. Similarly, in order for the 
 component which maps \(B_3\) to \(A_1\) to be nontrivial, it must send
 \(1 \in B_3 \cong \Q[X]/(X^{N-1})\) to something nonzero. Again, a
 consideration of the \(\Delta\) grading shows this is impossible. 

Finally, if both \(L_1\) and \(L_3\) are links, we must consider the
component of \(f_3\) which maps \(B_3\) to \(B_1\). Since \(f_3\) is a
map of \(\Q[X]\) modules, this homomorphism must be equal to multiplication by
some polynomial \(p(X)\). Inspecting the bigrading on the two
summands we find that we must have \(p(X) = cX\).
%
%
If \(c=0\), then \(B_1\) injects into \(\hkrn(L_2)\), so 
the action of \(X\) on this group is nontrivial. But this is
impossible, since \(L_2\) is a knot.  
\end{proof}

Now that we understand the action  of \(f_3\), it is
straightforward to determine  \(\khrn(L_2)\)
from \(\khrn(L_1)\) and \(\khrn(L_3)\) and to check that it has the
expected form. We give a detailed argument in the case where
\(L_1 = \undercrossing\), \(L_2 =
\asmoothing\), and \(L_3 = \doublepoint\), 
and leave the other cases (which are similar) to the reader.

Suppose first that \(L_1 \) is a knot and 
 \(L_2 \) is a regular link. Then from equation~\eqref{Eq:SingHOM}, we see that
\begin{align}
P(L_2) & = a q^{-1}P(L_1) + q^{-1}P(L_3)  \\
       & = a q^{-1}P(L_1) + q^{-1}\wP(L_3) + q^{-1}Q(L_3) 
\end{align}
Since \(\lk L_2 = \lk L_3\) and \(\sigma (L_2) = \sigma (L_3) -1\), the
term \( q^{-1}Q(L_3)\) is equal to \(Q(L_2)\), which means that 
\begin{equation}
\label{Eq:L2Pt}
\wP(L_2) =  a q^{-1}P(L_1) + q^{-1}\wP(L_3). 
\end{equation}

On the other hand, the fact that \(L_2\) is a link implies that \(f_3
= 0\), so the skein exact sequence splits to give a short exact sequence
\begin{equation*}
\begin{CD}
0 @>>> \hkrn(L_1) @>{(0,N-1)}>> \hkrn(L_2) @>{(0,1)}>>
\hkrn(L_3) @>>> 0
\end{CD}
\end{equation*}
from which we get the corresponding equation
\begin{align}
\khrn(L_2) & = q^{N-1} \khrn(L_1) + q^{-1}\khrn(L_3) \\
           & =  (aq^{-1} \mathscr{P}(L_1) + 
q^{-1} \wcP(L_3))\vert_{a=q^N} + q^{-1}\wQ(L_3) 
\end{align}
All the terms in \(aq^{-1} \mathscr{P}(L_1)\) have \(\delta\)-grading
\(1+ \sigma(L_1) = \sigma(L_2)\). Similarly, all the terms in \(
q^{-1} \wcP(L_3)\) have \(\delta\)-grading \(-1+\sigma(L_3) =
\sigma(L_2)\). Combined with equation~\eqref{Eq:L2Pt}, this implies
that 
\begin{equation}
\wcP(L_2) = aq^{-1} \mathscr{P}(L_1) + 
q^{-1} \wcP(L_3)
\end{equation}
so
\begin{equation}
\khrn(L_2) = \wcP(L_2) \vert_{a=q^N} + \wQ(L_2).
\end{equation}
This implies both that \(\khrn(L_2)\) has the expected form and that
the Laurent polynomial \(\wP(L_2)\) is alternating.

The case when \(L_2\) is a knot is somewhat more interesting. We have
\begin{align}
P(L_2) & = a q^{-1}P(L_1) + q^{-1}P(L_3)  \\
       & = a q^{-1}\wP(L_1) + a q^{-1}Q(L_1)  + q^{-1}\wP(L_3) + q^{-1}Q(L_3). 
\end{align}
Using the identities
 \( \lk L_1 + 1/2 =  \lk L_3\), \(\sigma(L_1) =
\sigma(L_2)-1\), and  \(\sigma(L_3) = \sigma(L_2) + 1\), we see that
\begin{align*}
a q^{-1}Q(L_1) + q^{-1}Q(L_3) & = 
\biggl[(aq^{-1})q^{\sigma(L_1)}(-a)^{2 \lk L_1} + 
q^{-1}q^{\sigma(L_3)}(-a)^{2 \lk  L_3} \biggr]
\left( \frac{aq^{-1}
  -a^{-1}q}{q-q^{-1}} \right) \\
& = (-q^{-2} + 1) q^{\sigma(L_2)}(-a)^{2 \lk L_3}
\left( \frac{aq^{-1}
  -a^{-1}q}{q-q^{-1}} \right) \\
& = (aq^{-2} - a^{-1})q^{\sigma(L_2)}(-a)^{2 \lk L_3}. 
\end{align*}
Thus 
\begin{equation}
P(L_2) = aq^{-1}\wP(L_1) + q^{-1}\wP(L_3) + (aq^{-2} - a^{-1})q^{\sigma(L_2)}(-a)^{2 \lk L_3}.
\end{equation}

The corresponding statement on the level of homology can be derived
from the short exact sequence
\begin{equation*}
\begin{CD}
0 @>>> A_1 \oplus B_1/XB_1@>{(0,N-1)}>> \hkrn(L_2) @>{(0,1)}>>
A_3 \oplus X^{N-2}B_3@>>> 0.
\end{CD}
\end{equation*}
We get
\begin{equation*}
\begin{split}
\khrn(L_2) & = \bigl[aq^{-1} \wcP(L_1) + q^{-1}\wcP(L_3) \\ & \quad + (aq^{-1})
q^{\sigma(L_1)}(at^{-1})^{2 \lk L_1} (a^{-1}q^{2}) +
(q^{-1})q^{\sigma(L_3)} (at^{-1})^{2 \lk L_3} (aq^{-2})\bigr ] \bigl
\vert_{a=q^N}  \\
& = \bigl[ aq^{-1} \wcP(L_1) + q^{-1}\wcP(L_3) +
(a^{-1}t+aq^{-2})q^{\sigma(L_2)}(at^{-1})^{2\lk L_3}\bigr ] \bigl
\vert_{a=q^N} \\
& = \mathscr{P}(L_2)\vert_{a=q^N}
\end{split}
\end{equation*}
(Again, it is easy to check that all the terms in the next-to-last
line have \(\delta\)-grading \(\sigma(L_2)\).) 

To complete the proof of the theorem,
it remains to check that the \(\Q[X]\) module
structure on \(\hkrn(L_2)\) agrees with that of a thin link. This is true if
\(L_2\) is a knot, since \(X\) always acts trivially in this case.
 If \(L_2\) is a
link, then exactly one of \(\hkrn(L_1)\) and \(\hkrn(L_3)\) contains a
\(\Q[X]/(X^{N-1})\) summand, so \(\khrn(L_2)\) has a sub- or
quotient module \(B \cong \Q[X]/(X^{N-1})\). For the sake of argument,
suppose \(B\) is a submodule. The \(\Delta\)-grading of \(1 \in B\) is
\(0\), so if \(B\) were contained in a direct
summand larger than itself, \(\hkrn(L_2)/B\) would contain an element
of \(\Delta\)-grading \(-2\). But every element of \(\hkrn(L_2)/B\)
has \(\Delta\)-grading \(0\), so for \(N>4\), this is a contradiction.
Thus \(B\) is a direct summand. 
Likewise, if the action of \(X\) on \(\hkrn(L_2)/B\) was nontrivial,
it would contain an element of \(\Delta\)-grading \(2\).
 This proves the claim about the module structure of
\(\hkrn(L_2)\)  when \(B\) is a submodule. The case of a quotient
module is similar. 
\end{proof}

\subsection{Twisting}

\begin{figure}
\includegraphics{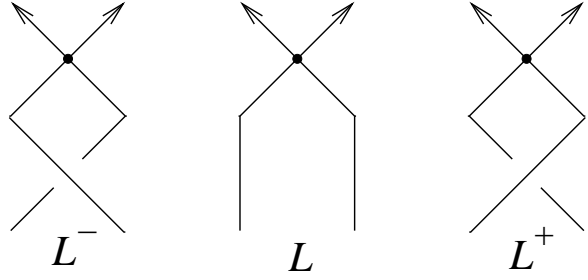}
\caption{\label{Fig:Untwist} Adding a twist to a singular link.}
\end{figure}
If \(L\) is a singular link, then we can modify \(L\) to produce new
links \(L^+\) and \(L^-\) by adding a positive or negative twist
adjacent to the singular point, as shown in Figure~\ref{Fig:Untwist}. 
(Warning: \(L^{\pm}\) should not be confused with \(L_{\pm}\) --- the
first is a singular link, while the second is a regular one.)
The twist can be added either below the singular point, as shown in
the figure, or above it (reverse all the orientations.) If two
singular links can be related by a sequence of such operations, we say
they are {\it twist equivalent}. 

 Using the second Reidemeister move, it is easy to
see that \((L^+)^- = L\). Thus if we want to study the effect of
twisting on an invariant of singular links, 
it suffices to consider the operation
of replacing \(L\) by \(L^-\). 
For example, it is not difficult to compute that 
\begin{align}
P(L^-) & = q^{-1} P(L_-) - a^{-1} P(L_0) \\
& = -a^{-1}q^{-1}(qP(L_0) - aP(L_-)) \\
\label{Eq:Twist}
& = -a^{-1}q^{-1}P(L). 
\end{align}
Its effects on \(\hkrn\) are similarly
mild.
\begin{lem}
\label{Lem:Untwist}
\( H_N(L_-) \cong H_N(L){[ 1 ] \{-N-1\}}\).
\end{lem}

\noindent{\bf Remark:} As in \cite{KRI}, the terms in  
brackets and braces indicate
shifts in the \(i\) and \(j\) gradings respectively, so if \(x \in
H_N(L)\) and \(x'\) is the corresponding element in \(H_N(L_-){[ 1 ] \{-N-1\}}\), we
have \(i(x') = i(x) + 1\), \(j(x') = j(x) - N-1\). 

The proof of the lemma is essentially contained in the proof that
\(H_N\) is invariant under the second Reidemeister move given in
\cite{KRI}. 

\begin{figure}
\input{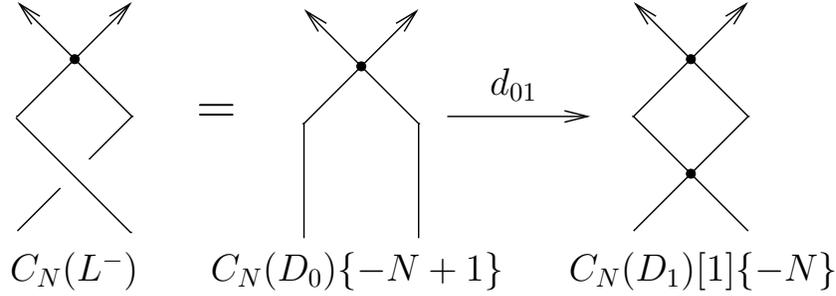}
\caption{\label{Fig:ReidII} The chain complex \(C_N(L^-)\).}
\end{figure}

\begin{proof}
The chain complex \(C_N(L^-)\) is shown schematically in
Figure~\ref{Fig:ReidII}. By Proposition 30 of \cite{KRI}, we know that
the chain complex 
\(C_N(D_1)\) can be decomposed  as \(C_N(D_0){\{1\}} \oplus
  C_N(D_0){\{-1\}}\). Let \((\alpha,\beta)\) be the components of 
\begin{equation}
d_{01}: C_N(D_0){\{-N+1\}} \to C_N(D_1){[1]\{-N\}} \cong C_N(D_0){[1]
 \{-N+1\}} \oplus
C_N(D_0){[1]\{-N-1\}}
\end{equation}
  By Lemma 25 of \cite{KRI},
  \(\alpha\) is an isomorphism. The claim now follows from a standard
  cancellation argument. Explicitly, we
  observe that \(C_N(L^-)\) has an acyclic subcomplex of the form
  \((C_N(D_0){\{-N+1\}}, \text{im} \ts d_{01})\), and that the quotient complex 
 \begin{equation}
\frac{C_N(D_1){[ 1 ]\{-N\}}}{  \text{im} \ts d_{01}} \cong
\frac{C_N(D_0){[ 1 ]\{-N+1\}} \oplus C_N(D_0){[ 1
  ]\{-N-1\}}}{(\alpha(x), \beta(x))}
\end{equation}
 is isomorphic to 
\(C_N(D_0){[ 1 ]\{-N-1\}}\) via the map which sends a pair
\((y,z)\) to \(z-\beta \alpha^{-1}(y).\) 
\end{proof}

Next, we investigate the effect of adding a twist on the linking
number and signature of \(L\). 

\begin{lem}
 \(\lk  L^- = \lk L - 1/2\) and \(\sigma(L^-) =
\sigma(L) - 1\). 
\end{lem}

\begin{proof} The first equation is elementary. To
  prove the second, we substitute \(a=-1\),
  \(q=i\) into equations~\eqref{Eq:SingHOM} and \eqref{Eq:Twist} to obtain
\begin{align}
\Det L \phantom{^-} & = i \Det L_0 + \Det L_- \\
\Det L^- & = -i \Det L.
\end{align}
The problem now breaks into three cases, depending on the relative
values of \(\Det L_-\) and \(\Det L_0\). We consider each case
separately.
\vskip0.05in
\noindent{\bf Case 1:} \(\phi(L_-) = i \phi(L_0)\).
 Then \(\phi(L) = i \phi(L_0) = \phi(L_-)\),  so 
\(\sigma(L) = \sigma(L_0)+1 = \sigma(L_-)\). 
Also \(\phi(L^-) = - i \phi(L) = -i 
\phi(L_-)\), so \(\sigma(L^-) = \sigma(L_-)-1 = \sigma(L) -1\). 
\vskip0.05in
\noindent{\bf Case 2:} \(\phi(L_-) = -i \phi(L_0)\) and \(\det L_- >
\det L_0\). Then \(\phi(L) = -i \phi(L_0) = \phi(L_-)\), so
\(\sigma(L) = \sigma(L_0) - 1 = \sigma(L_-)\). Also \(\phi(L^-) = -i
\phi(L) = - i \phi(L_-)\), so \(\sigma(L^-) = \sigma(L_-) -1 =
\sigma(L) -1\).
\vskip0.05in
\noindent{\bf Case 3:} \(\phi(L_-) = -i \phi(L_0)\) and \(\det L_- <
\det L_0\). Then \(\phi(L) = i \phi(L_0) = - \phi(L_-)\), so \(\sigma(L) =
\sigma(L_0)+1\), \(\sigma(L_-) = \sigma(L_0)-1\). We have \(\phi(L^-) =
-i \phi(L) = i \phi(L_-)\), so \(\sigma(L^-) = \sigma(L_-)+1 = \sigma(L)-1.\)
\vskip0.05in
\noindent{\bf Case 4:} \(\det L_- = 0 \). Then \(\phi(L) = i
\phi(L_0)\), so \(\sigma(L) = \sigma(L_0) + 1\). On the other hand, we
have \(\sigma (L_-) = \sigma (L_0)\) and \(\phi(L^-) = -i \phi(L) =
\phi(L_0)\), so \(\sigma(L^-) = \sigma(L_-) = \sigma(L_0)\). 
\vskip0.05in
\noindent{\bf Case 5:} \(\det L_0 = 0 \). Then \(\sigma(L) =
\sigma(L_0) = \sigma(L_-)\) and \(\phi(L^{-}) = -i \phi(L) = -i
\phi(L_-)\), so \(\sigma(L^-) = \sigma(L_-) - 1.\)

\end{proof}

\noindent Putting these facts together, we obtain 
\begin{cor}
\label{Cor:UnTwist}
Let \(L_1\) and \(L_2\) be twist-equivalent singular links.
 Then \(L_1\) is \(N\)-thin if and only
if \(L_2\) is. 
\end{cor}

\section{Some Examples}
\label{Sec:Examples}

We conclude by using the results of the previous section to show  that many
small knots and links are \(N\)-thin for \(N > 4\).
 Two-bridge knots and links provide the
best class of examples, but the method can also be applied to 
other knots with small crossing number. In section~\ref{Sec:Other}, we
give  examples  of such knots, as well as a few interesting
knots which are not \(N\)-thin. Finally, in
section~\ref{SubSec:Unreduced}, we combine our results on the reduced
homology with a theorem of Gornik to compute the unreduced
Khovanov-Rozansky homology in a few cases.

\subsection{Two-bridge links}
\begin{figure}
\includegraphics{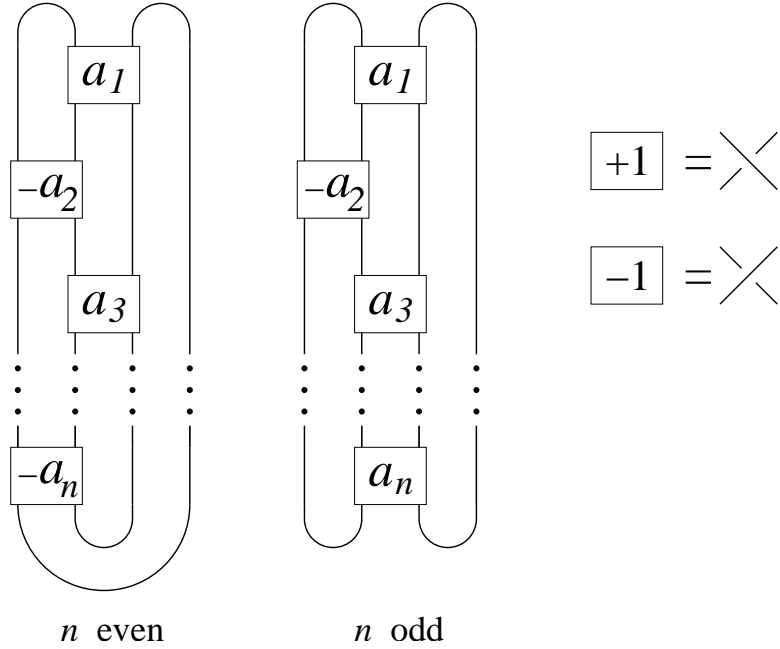}
\caption{\label{Fig:TwoBridge} Plat diagrams for two-bridge links.}
\end{figure}
\label{Sec:TwoBridge}
We are now in position to complete the proof of
Theorem~\ref{Thm:One}. 
We begin by recalling a few basic facts about two-bridge knots and
links. This is all standard material --- see {\it e.g.}
\cite{Murasugi} for a more detailed exposition. 

 A link in \(S^3\) is two-bridge if and only if it
 has a planar diagram of the form shown in
 Figure~\ref{Fig:TwoBridge}.
Such links are classified up to isotopy by pairs of relatively prime integers
\((p,q) \ (p\geq 0)\) modulo the relation \((p,q)\sim (p,q')\) if \(q
\equiv q' \pmod p\) or \(q' \equiv q^{-1} \pmod p\). 
The integers \(a_i\) are the coefficients
 in a continued fraction expansion for \(p/q\):
\begin{equation}
\label{Eq:Cfrac}
\frac{p}{q} = a_1 + \cfrac{1}{a_2+\cfrac{1}{a_3+\cfrac{1}{\ldots +
      \cfrac{1}{a_n}}}}
\end{equation}
(If \(n=0\), the link is the unlink, and  \((p,q)\) is defined to
be \((0,0)\).)
We write \(K(p,q)\)
to denote the two-bridge link determined by the pair
\((p,q)\). If \(p\) is odd, \(K(p,q)\) is a knot, while if \(p\) is
even, it is a two-component link, both of whose components are
unknots. Finally, \( \det K(p,q) = p. \)

It is easy to see that the diagram of Figure~\ref{Fig:TwoBridge} is
alternating if and only if all the \(a_i\)'s have the same sign. If
\(0 <q < p\), it is well known that \(p/q\) admits a continued
fraction expansion of the type shown in equation~\eqref{Eq:Cfrac} with
all \(a_i >0\). Thus any two-bridge knot or link can be represented by an
alternating two-bridge diagram. 

\begin{lem}
\label{Lem:TBRes}
 Let \(D\) be an alternating two-bridge diagram  representing
\(K(p,q)\)
  \((0 < q < p)\). Then the links obtained by resolving the top
  crossing of \(D\) in both possible ways are 
 \(K(q,p)\) and  \(K(p-q,q)\). 
\end{lem}

\begin{proof}

\begin{figure}
\includegraphics{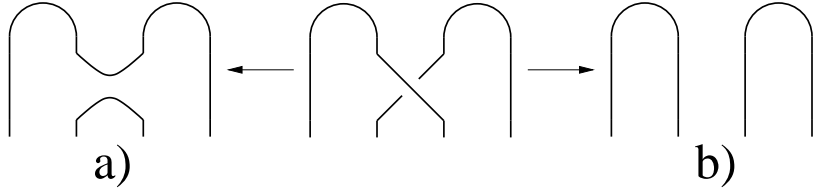}
\caption{\label{Fig:TwoRes} Resolving the top crossing in a two-bridge
  link.}
\end{figure}

By hypothesis, \(D\) is an alternating diagram, so all the \(a_i\)'s
have the same sign. Without loss of generality, we assume that they
are all positive.
If we resolve the top crossing as in 
 Figure~\ref{Fig:TwoRes}a, then we can untwist the remaining
\(a_1-1\) twists in the central column,
 move the rightmost strand all the way to the left,
and flip the whole diagram over to obtain a two-bridge diagram
\(K(p_1,q_1)\) with
continued fraction expansion \([a_2,a_3, \ldots,a_n]\). We have
\begin{equation}
\frac{p}{q} = a_1 + \frac{q_1}{p_1}
\end{equation}
so \((p_1,q_1) = (q, p \pmod q)\).

On the other hand, if the top crossing is resolved as 
 in Figure~\ref{Fig:TwoRes}b
and \(a_1>1\), then we obviously have a two-bridge diagram
 \(K(p',q')\) with
continued fraction expansion \([a_1-1,a_2,\ldots,a_n]\). Thus
\begin{equation}
\frac{p}{q} = a_1 + \frac{q_1}{p_1} \qquad \frac{p'}{q'} = a_1 -1 + 
\frac{q_1}{p_1}
\end{equation}
so \((p',q') = (p-q,q)\). 

Finally, if the resolution is as in
Figure~\ref{Fig:TwoRes}b and \(a_1 = 1\), then we undo the \(a_2\)
twists in the top left of the resolved diagram and are left with a
two-bridge diagram \(K(p_2,q_2)\)
with continued fraction expansion \([a_3,a_4,
  \ldots,a_n]\). Here we have
\begin{equation}
\frac{p}{q}  = 1 + \frac{1}{a_2+{q_2}/{p_2}} 
 = \frac{(a_2+1)p_2 + q_2}{a_2 p_2 + q_2}.
\end{equation}
so  \((p_2,q_2) = (p-q, q \pmod {p-q})\).

\end{proof}

\begin{defn}
A {\it special singular two-bridge link} is a singular link obtained
by replacing the top crossing of an alternating two-bridge diagram
with a singular crossing. 
\end{defn}

\begin{thrm}
Suppose \(L\) is either a two-bridge knot or link or a special singular
two-bridge link  and that \(\det L \neq 0
\).  Then \(L\) is \(N\)-thin (with respect to either component)
for all \(N>4\). 
\end{thrm}

\noindent {\bf Remark:} The condition \(\det L \neq 0\) rules out the
unlink, which we have already seen is not thin.

\begin{proof}
By induction on \(\det L\). For the base case, we must consider links
and singular links with \(\det L = 1\). There is a unique regular
two-bridge link with determinant \(1\), namely the unknot, so
the claim is easily verified in this case. 
Now suppose \(L\) is a singular link. The unoriented resolution
\(L_u\) is a two-bridge diagram of a knot with
determinant \(1\), so it must be the unknot. Moreover \(L_u\) is
alternating, and  it is well
known that any alternating diagram of the unknot is {\it nugatory} ---
it can be reduced to the standard diagram by repeated application of
the Reideimester I move. Thus \(L\) is twist-equivalent to 
the singular link \(\Theta\) illustrated in Figure~\ref{Fig:Theta}. 
By Corollary~\ref{Cor:UnTwist}, it suffices to show that \(\Theta\) is
\(N\)-thin. 

\(\Theta\) is clearly symmetric, so without loss of
generality, we assume that the marked component belongs to the side
labeled \(X_1\). By Lemma~\ref{Lem:Theta}, we have
\begin{equation}
\hkrn(\Theta,C_1) \cong \A_N(\Theta)/(X_1)\cong \Q[X_2]/(X_2^{N-1}).
\end{equation}
and the Poincar{\'e} polynomial is
\begin{equation}
\khrn(\Theta,C_1) = q^{-N+2} + q^{-N+4} + \ldots + q^{N-4} + q^{N-2}.
\end{equation}
On the other hand, 
\begin{equation}
P(\Theta) = q\frac{a-a^{-1}}{q-q^{-1}} - a = \frac{aq^{-1}-a^{-1}q}{q-q^{-1}}
\end{equation}
so \(\Det \Theta = 1\). Since \(\Theta_0\) is the unlink, which has
both determinant and signature equal to \(0\), we see that
\(\sigma(\Theta)=\lk(\Theta) = \wP(\Theta)=0\). 
It follows that \(\Theta\) is \(N\)-thin.

We now assume the theorem holds for links with determinant \(<p\).  
Given a regular two-bridge link \(L=K(p,q) \ (0 < q < p)\),
 we choose an alternating two-bridge diagram
representing \(L\), and let \(L_0\) and \(L_s\) denote the links
obtained by taking the oriented and singular resolutions of its
topmost crossing, so that \(L\), \(L_0\), and \(L_s\) fit into a skein
exact sequence. 
 \(L_s\) is clearly a special singular two-bridge link, and  \(L_0\) is a
two-bridge link by Lemma~\ref{Lem:TBRes}.
Combining 
Lemmas~\ref{Lem:DetRes} and \ref{Lem:TBRes}, we see that one of
\(L_0,L_s\) has determinant \(q\) and the other has determinant
\(p-q\). 
 Thus both \(L_0\) and \(L_s\) are \(N\)-thin by the induction
hypothesis, and 
\begin{equation}
 \det L =  \det L_0 + \det L_s. 
\end{equation}
Then Theorem~\ref{Thm:ThinSeq}
 implies that \(L\) is \(N\)-thin as well.

\begin{figure}
\includegraphics{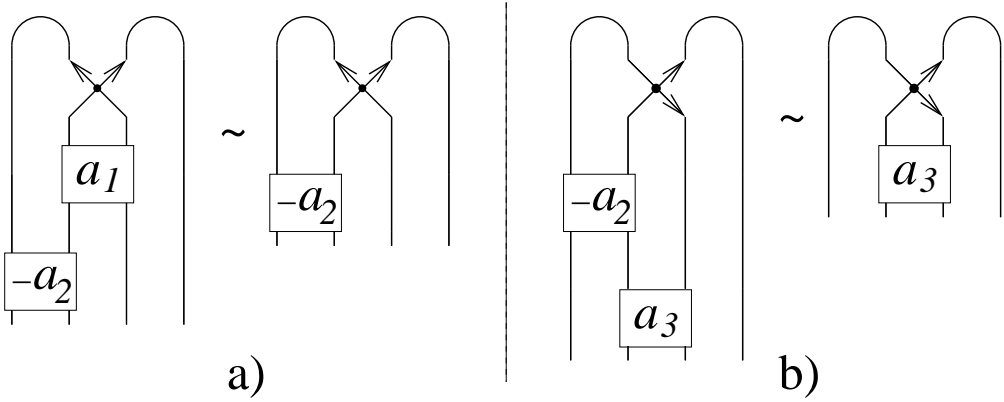}
\caption{\label{Fig:SingOr} Possible orientations for the singular crossing.}
\end{figure}

Next, suppose we are given a special
singular two-bridge link \(L\) with determinant
\(p\). We distinguish two cases, depending on the orientation of the
singular crossing. First, suppose that the crossing is oriented as in
Figure~\ref{Fig:SingOr}a or its reverse.  Then as shown in the
figure, \(L\) is 
twist-equivalent to a diagram with just the singular crossing in the
top block. The unoriented resolution \(L_u\) is a two-bridge link
\(K(p,q)\) with continued fraction expansion \([a_2,a_3,\ldots,a_n]\),
where all the \(a_i\)'s have the same sign. Suppose for the moment
that \(a_i>0\), so \(0 < q < p\). There is a skein exact sequence
relating \(L\), \(L_0 = K(p_0,q_0)\), and \(L_+ = K(p_+,q_+)\), where
 the  continued fraction expansions of \(p_0/q_0\) and \(p_+/q_+\) are
 \([a_3,a_4,\ldots,a_n]\) and \([-1,a_2,a_3,\ldots,a_n]\). Thus 
\begin{equation}
\frac{p}{q}  = a_2 + \frac{q_0}{p_0} 
 = \frac{a_2p_0+q_0}{p_0}
\end{equation}
and
\begin{equation}
\frac{p_+}{q_+}  = -1 + \frac{1}{a_2 +{q_0}/{p_0}} 
 = \frac{(a_2-1)p_0 + q_0}{a_2p_0+q_0}.
\end{equation}
In other words, 
\begin{align}
\det L\phantom{_0} & =  a_2p_0 + q_0 \\
\det L_0 & = p_0 \\
\det L_+ & = (a_2-1)p_0 + q_0
\end{align}
so \(\det L = \det L_0 + \det L_+\). The induction hypothesis implies
that both \(L_0\) and \(L_+\) are \(N\)-thin, so by
Theorem~\ref{Thm:ThinSeq}, we conclude that \(L\) is \(N\)-thin as
well. The argument when \(a_i <0\) is identical, except that we
replace \(L_+\) with \(L_-\). 

We now consider the case where the singular crossing is oriented
as in Figure~\ref{Fig:SingOr}b or its reverse. We claim that without
loss of generality, we may assume \(a_1 \neq 0 \). Indeed, suppose this is not
the case. Then as indicated
in the figure, \(L\) is twist-equivalent to a singular two-bridge knot
with continued fraction expansion \([a_3,a_4,\ldots,a_n]\). As in
the previous case, all the \(a_i\)'s have the same sign; we assume for
the moment that \(a_i>0\). 

In this case, we consider the skein exact sequence relating \(L\),
 \(L_-\), and \(L_0\). Let  \(L_0 = K(p_1,q_1)\), where
 \(p_1/q_1\) has continued fraction expansion \([a_2,a_3,\ldots,
 a_n]\). The 
unoriented resolution \(L_u\) has continued fraction expansion
\([a_1,a_2,\ldots, a_n]\), so \(\det L  = a_1p_1 + q_1\). On
the other hand, \(L_-\) has continued fraction expansion
\([a_1-1,a_2, \ldots,a_n]\). Since \(a_1>0\), it follows that
 \(\det L_+ = (a_1-1)p_1 + q_1\). Thus \(\det L = \det L_0 + \det
 L_+\), and \(L\) is \(N\)-thin. 
Finally, when \(a_i<0\), we argue as above, but with \(L_-\) replaced
by \(L_+\). This completes the proof of the theorem. 

\end{proof}
\subsection{Other Knots}
\label{Sec:Other}
Theorem~\ref{Thm:ThinSeq} and Corollary~\ref{Cor:UnTwist} can also be
applied on a case-by-case basis to show that certain other 
knots are thin. For example, we have the following 
\begin{crit}
\label{Crit:Change}
Suppose \(L_1\) and \(L_2\) are two knots related by a crossing
change, and that 
\begin{align}
\phi(L_1) & = \phi(L_2) \\
\det L_1 & < \det L_2.
\end{align}
Then if \(L_1\) and the oriented resolution \(L_0\) are \(N\)-thin \((N>4)\),
\(L_2\) is \(N\)-thin as well. 
\end{crit}
\begin{figure}
\includegraphics{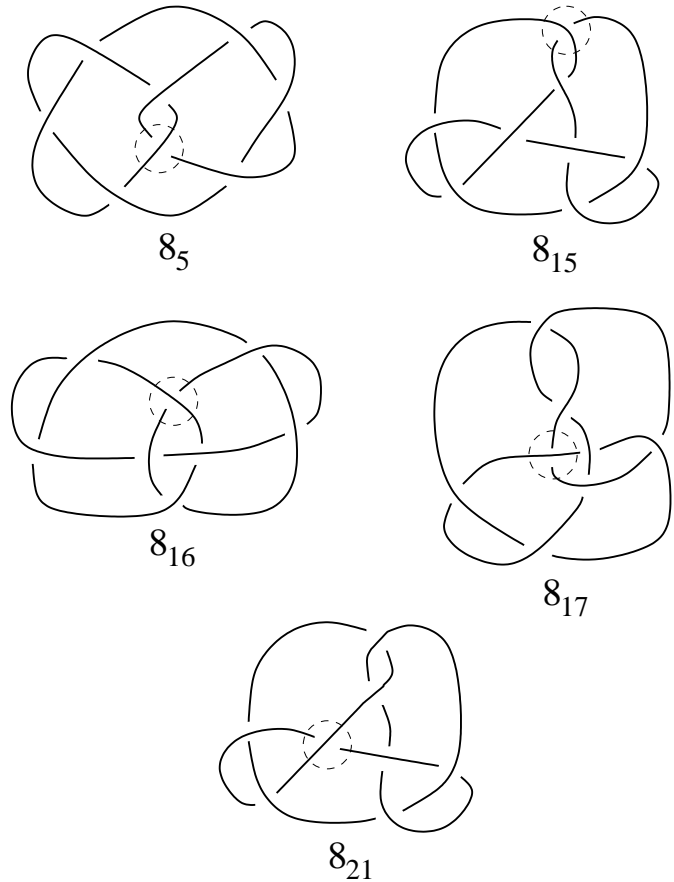}
\caption{\label{Fig:ThinKnots} \(8\)-crossing knots to which
  criterion~\ref{Crit:Change} can be applied. The relevant crossings
  are marked by circles. }
\end{figure}
\begin{proof}
We have 
\begin{align}
\Det L_- & = \Det L_s - i \Det L_0 \\
\Det L_+ & = \Det L_s + i \Det L_0 
\end{align}
so if \(L_-\) and \(L_+\) have the same phase, \( \det L_s > \det
L_0\). It follows that
\begin{align}
\det L_s & =  \det L_1 + \det L_0 \\
\det L_2 & = \det L_s + \det L_0
\end{align}
so two applications of Theorem~\ref{Thm:ThinSeq} give the desired
result. 
\end{proof}
This criterion provides a fast and reasonably effective way of finding
thin knots with small crossing number. 
Figure~\ref{Fig:ThinKnots} shows diagrams of five
\(8\)-crossing knots to which the criterion can be applied.
 Combining  with the results of
 Theorem~\ref{Thm:One}, we see that the 
  only knots with 8 crossings or fewer which are
not known to be \(N\)-thin (\(N>4\)) are 
\(8_{10},8_{18},8_{19}\), and \(8_{20}\). 
Of these, the knot \(8_{19}\) (the \((3,4)\) torus knot) cannot
 be \(N\)-thin for any \(N>2\), since its HOMFLY polynomial is not
 alternating. (Of course, it is not thin for \(N=2\) either.)
 It seems quite plausible that the remaining three knots
 are thin, although we do not have any proof of this fact. 

Interestingly, there also exist alternating knots whose HOMFLY
polynomials are nonalternating, although they are somewhat harder to
find. The smallest such knot is \(11^a_{263}\) ({\it Knotscape}
numbering). It has HOMFLY polynomial
\begin{multline*}
P(11^a_{263}) = a^{-8}(q^{-8}-q^{-6}+4 q^{-4} - 3q^{-2} + 6 -3q^2 +
4q^4-q^6+q^8) \\
+ a^{-10}(q^{-8}-4q^{-6}+4q^{-4}-9q^{-2}+5-9q^2+4q^4-4q^6+q^8) \\
+ a^{-12}(-q^{-6}+3q^{-4}-2q^{-2} + 5 -2 q^2 + 3 q^4 - q^6) 
- a^{-14}
\end{multline*}
which is alternating except at the final term. This knot cannot be
\(N\)-thin for any \(N>2\), but
its ordinary Khovanov homology is thin by Lee's theorem. 

\subsection{Unreduced Homology} 
\label{SubSec:Unreduced}
So far, we have been concerned only with the reduced version of the
Khovanov-Rozansky homology. We conclude by briefly considering
to the problem of computing the unreduced homology. Here,
unfortunately, our methods are  less successful. Nonetheless,
combining our knowledge of the reduced homology with a theorem of
Gornik \cite{Gornik} enables us to compute the unreduced homology in a
few special cases. The first ingredient in the calculation is a spectral
sequence relating the two theories. 

\begin{lem}
\label{Lem:RedtoUnred}
There is a spectral sequence with \(E_1\) term \(\hkrn(L,C)\otimes
\Q[X]/(X^N)\) which converges to \(H_N(L)\). All of its differentials
respect the \(q\)-grading. 
\end{lem}

\begin{proof}
Let \(X\) be the edge operator corresponding to the component \(C\). 
The action of \(X\) on \(C_N(L)\) induces a filtration 
\(C_N(L) = F_0 \supset F_1 \supset F_2 \ldots \supset F_N =\{0\}\)
where \(F_i = X^i C_N(L)\). From property (4) of \ref{Prop:Foams} it
follows that \(F_i/F_{i+1} \cong C_N(L)/XC_N(L)\) for \(0\leq i <
N\). The homology of \( C_N(L)/XC_N(L)\) is \(\khrn(L,C)\), so the
\(E_1\) term of the associated spectral sequence is 
\(\hkrn(L,C)\otimes \Q[X]/(X^N)\). Finally
the fact that the differentials in this spectral
sequence respect the \(q\)-grading merely reflects the fact that the
same is true in \(C_N(L)\). 
\end{proof}

  We will also need Gornik's generalization of Lee's spectral sequence to the
Khovanov-Rozansky homology~\cite{Gornik}. 
\begin{thrm} (Gornik)
Let \(K\) be a knot.  
There is a spectral sequence with \(E^1\) term \(H_N(K)\) which
converges to a vector space of dimension \(N\) supported in
homological grading \(0\). The differential \(d_i\) lowers the
\(q\)-grading by \(2iN\). 
\end{thrm}

\begin{prop}
Let \(T_{2,n}\) be the positive \((2,n)\) torus knot. For \(N>4\),
 the Poincar{\'e}
polynomial of \(H_N(T_{2,n})\) is given by
\begin{equation}
\mathscr{P}(H_N((T_{2,n})) = q^{(n-1)(N-1)}\left[
  [N]+[N-1]q^{-1}(1+q^{2N}t^{-1}) 
\sum_{i=1}^{(n-1)/2} q^{4i}t^{-2i}   \right] 
\end{equation}
where as usual
\begin{equation}
[N] = \frac{q^{N}-q^{-N}}{q-q^{-1}}.
\end{equation}
\end{prop}

\begin{proof}
The HOMFLY polynomial of \(T_{2,n}\) is given by
\begin{equation}
P(T_{2,n}) = (aq^{-1})^{n-1} \left[ \sum_{i=0}^{(n-1)/2} q^{4i} - 
a^2q^2 \sum_{i=0}^{(n-3)/2} q^{4i}  \right]
\end{equation}
and \(\sigma(T_{2,n}) = n-1\), so 
\begin{equation}
\khrn(T_{2,n}) = (aq^{-1})^{n-1}\left[ \sum_{i=0}^{(n-1)/2} q^{4i}t^{-2i} + 
a^2q^2t^{-3} \sum_{i=0}^{(n-3)/2} q^{4i}t^{-2i}  \right]
_{a=q^N} 
\end{equation}
In other words, \(\hkrn(T_{2,n})\) is generated by classes \(a_i\)
(\(0\leq i \leq (n-1)/2\)) and \(b_j\) (\(0\leq j \leq (n-3)/2\)),
with homological gradings \(-2i\) and \(-2j-3\) respectively. In
particular, there is one generator in grading \(0\) and in
each of the gradings \(-2,-3,\ldots,-n\). 

Consider the generators \(a_1
\otimes X^k\) in the \(E_1\) term of the spectral sequence of
Lemma~\ref{Lem:RedtoUnred}. They have homological grading
\(-2\), so they can only be killed by generators of the form \(b_0
\otimes X^l\). Since the differentials in this
spectral sequence respect the \(q\)-grading, it is not difficult to
check  that the only
possible differential from \(b_0 \otimes X^l \) to \(a_1 \otimes X^k\)
takes \(b_0 \otimes 1\) to \(a_1 \otimes X^{N-1}\). We claim that this
differential must be nontrivial. Indeed, suppose it were not. Then the
\(E_1\) term of Gornik's spectral sequence would contain the generator
\(a_1 \otimes X^{N-1}\), which must be killed by some generator with
homological grading \(-3\), and the \(q\) grading of this element would be
at least \(2N\) greater than the \(q\)-grading of \(a_1 \otimes
X^{N-1}\). But the generator in grading \(-3\) with largest
\(q\)-grading is \(b_0 \otimes X^{N-1}\), and its \(q\)-grading is
only \(2N-2\) greater than that of \(a_1 \otimes X^{N-1}\). Thus we
have arrived at a contradiction. 

It follows that in homological grading \(-2\), 
\(H_N(T_{2,n})\) has generators \(a_1 \otimes X^k\)
(\(0\leq k \leq N-2\)). In Gornik's spectral sequence, these must be
killed by the generators 
\(b_0 \otimes X^l\) (\(1\leq l \leq N-1\)). In particular, all
of these generators must survive in the reduced-unreduced spectral
sequence, so there are no differentials from \(a_2\otimes X^k\) to
\(b_0 \otimes X^l\).  
To finish the proof, we simply repeat this argument, considering
differentials from \(b_1 \otimes X^l\) to \(a_2 \otimes X^k\), then
from \(b_2 \otimes X^l\) to \(a_3 \otimes X^k\), and so on. 
\end{proof}

This result is in accordance with the behavior predicted in section 5.10 of
\cite{superpolynomial}. 
A similar argument can also be used to compute the unreduced homology
of the figure-eight knot. For \(N>4\), we find
\begin{equation}
\mathscr{P}(H_N(4_1)) = [N] + [N-1](q^{2N+1}t^{-2} + qt^{-1} + q^{-1}t + q^{-2N-1}t^2)
\end{equation}
thus validating the prediction made in equation (61) of
\cite{superpolynomial}.

\bibliographystyle{plain}
\bibliography{../mybib}

\end{document}